\documentclass[leqno]{article}
\usepackage{mfpic}
\usepackage{amsmath,amsthm,amssymb}
\usepackage[T1]{fontenc}
\opengraphsfile{ratrys}

\newtheorem{Th}{Theorem}
\newtheorem{Prop}[Th]{Proposition}
\newtheorem{Lemma}[Th]{Lemma}
\newtheorem{Cor}[Th]{Corollary}

\theoremstyle{remark}
\newtheorem{remark}{Remark}
\newtheorem{definition}{Definition}
\newtheorem*{proofpr}{Proof of Theorem~\ref{prpr}}
\newtheorem*{problem}{Problem}

\newcommand{\cs}{{\cal S}}
\newcommand{\ct}{{\cal T}}

\newcommand{\st}{(S_t)_{t\in\R}}
\newcommand{\ttt}{(T_t)_{t\in\R}}

\newcommand{\s}{\mathbb{S}}
\newcommand{\T}{{\mathbb{T}}}

\newcommand{\R}{{\mathbb{R}}}
\newcommand{\C}{{\mathbb{C}}}
\newcommand{\Z}{{\mathbb{Z}}}
\newcommand{\N}{{\mathbb{N}}}
\newcommand{\xbm}{(X,{\cal B},\mu)}
\newcommand{\ycn}{(Y,{\cal C},\nu)}

\newcommand{\vep}{\varepsilon}
\newcommand{\sgn}{\operatorname{sgn}}
\newcommand{\rat}{\operatorname{R}}
\newcommand{\var}{\operatorname{Var}}

\title{On mild mixing of special flows over irrational
rotations under piecewise smooth functions}
\author{Krzysztof Fr\k{a}czek \& Mariusz Lema\'nczyk}
\begin{document}
\maketitle
\renewcommand{\thefootnote}{}
\footnote{Research partly supported by KBN grant 1 P03A 038 26.}
\renewcommand{\thefootnote}{\fnsymbol{footnote}}

\begin{abstract} It is proved that all special flows
over the rotation by an irrational $\alpha$ with bounded partial
quotients and under $f$ which is piecewise absolutely continuous
with a non-zero sum of jumps are mildly mixing. Such flows are
also shown to enjoy a condition which emulates the Ratner
condition introduced in \cite{Rat}. As a consequence we construct
a smooth vector--field on $\T^2$ with one singularity point such
that the corresponding flow $(\varphi_t)_{t\in\R}$ preserves a
smooth measure, its set of ergodic components consists of a family
of periodic orbits and one component of positive measure on which
$(\varphi_t)_{t\in\R}$ is mildly mixing and is spectrally disjoint
from all mixing flows.
\end{abstract}

\section{Introduction}
The property of mild mixing of a (finite) measure--preserving
transformation has been introduced by Furstenberg and Weiss in
\cite{Fu-We}. By definition, a finite measure--preserving
transformation is mildly mixing if its Cartesian product with an
arbitrary ergodic (finite or infinite not of type I)
measure-preserving transformation remains ergodic. It is also
proved in \cite{Fu-We} that a probability measure--preserving
transformation $T:\xbm\to\xbm$ is mildly mixing iff $T$ has no
non-trivial rigid factor, i.e.\
$\liminf_{n\to+\infty}\mu(T^{-n}B\triangle B)>0$ for every
$B\in\mathcal{B}$, $0<\mu(B)<1$. For importance and naturality of
the notion of mild mixing see e.g.\
\cite{Aa-Li-We,Da-Le,Fu1,Le-Le,Le-Pa0,Sch-Wa}.

It is immediate from the definition that the (strong) mixing
property of an action implies its mild mixing which in turn
implies the weak mixing property. In case of Abelian non--compact
group actions, Schmidt in \cite{Sch} constructed examples (using
Gaussian processes) of mildly mixing actions that are not mixing.
A famous example of a mildly mixing but not  mixing system is the
well--known Chacon transformation $T$ (mild mixing of $T$ follows
directly from the minimal self-joining property of $T$,
\cite{Ju-Ra-Sw}). However, none of known examples of mild but not
mixing dynamical system was proved to be of smooth origin (see
\cite{Ka-Th} and \cite{Ha-Ka} - discussions about the three
paradigms of smooth ergodic theory).

Special flows built over an ergodic rotation on the circle and
under a piecewise $C^1$--function with non--zero sum of jumps were
introduced and studied by J.\ von Neumann in
\cite{Neu}\footnote[1]{We thank A.\ Katok for turning our
attention to this article.}. He proved that such flows are weakly
mixing for each irrational rotation. The weak mixing property was
then proved for the von Neumann class of functions but over
ergodic interval exchange transformations by Katok in
\cite{Ka-Ro}, while in \cite{Iw-Le-Ma} the weak mixing property
was shown for the von Neumann class of functions where the
$C^1$--condition is replaced by the absolute continuouity, however
in \cite{Iw-Le-Ma}, $T$ is again an arbitrary irrational rotation.
The absence of mixing of $T^f$ is well--known; it has been proved
by Ko\v{c}ergin in \cite{Ko}. In fact, from the spectral point of
view special flows in this paper  have no spectral measure which
is Rajchman, i.e.\ they are spectrally disjoint from mixing flows
(see \cite{Fr-Le}).

The aim of this paper is to show that the class of special flows
built from a piecewise absolutely continuous function $f:\T\to\R$
with a non--zero sum of jumps and over a rotation by $\alpha$ with
bounded partial quotients is mildly mixing (Theorem~\ref{twmmsf}).
One of the main tools, which yet can be considered as another
motivation of this paper, is Theorem~\ref{radlniez} in which we
prove that $T^f$ satisfies a property similar to the famous Ratner
property\footnote[2]{The possibility of having the Ratner property
for some special flows over irrational rotations was suggested to
us by B.\ Fayad and J.-P.\ Thouvenot.} introduced in \cite{Rat}
(see also \cite{Th}). It will follow  that any ergodic joining of
$T^f$ with any ergodic flow $(S_t)$ is either the product joining
or  a finite extension of $(S_t)$. In Section~\ref{partrig}, the
absence of partial rigidity for $T^f$ will be shown. Finally these
two properties combined will yield mild mixing (see
Lemma~\ref{jmm}).

As a consequence of our measure--theoretic results we will
construct a mildly mixing (but not mixing) $C^{\infty}$--flow
$(\varphi_t)_{t\in\R}$ whose corresponding vector--field has one
singular point (of a simple pole type). More precisely, we will
construct $(\varphi_t)_{t\in\R}$ on the two--dimensional torus,
such that $(\varphi_t)$ preserves a positive $C^{\infty}$--measure
and the family of ergodic components of $(\varphi_t)$ consists of
a family of periodic orbits and one non--trivial component of
positive measure which is mildly mixing but not mixing. More
precisely, the non--trivial component of $(\varphi_t)$ is
measure--theoretically isomorphic to a special flow $T^f$ which is
built over an irrational rotation $Tx=x+\alpha$ on the circle and
under a piecewise $C^{\infty}$--function $f:\T\to\R$ with a
non--zero sum of jumps. In these circumstances $T^f$ lies in the
parabolic paradigm (see \cite{Ha-Ka}).

Some minor changes in the construction of the $C^\infty$--flow
$(\varphi)_{t\in\R}$ (which uses some ideas descended from Blokhin
 \cite{Bl}) yield an ergodic $C^{\infty}$--flow which is mildly
mixing but not mixing and lives on the torus with attached
M\"{o}bius strip. This flow will enjoy the Ratner property in the
sense introduced in Section~\ref{ratprop}.

The authors would like to thank the referee for numerous remarks
and comments that improved the first version of the paper, and
especially for shortening the proof and for a strengthening of
Theorem~\ref{brakszt}.

\section{Basic definitions and notation}
Assume that $T$ is an ergodic automorphism of a standard
probability space $\xbm$. A measurable function $f:X\to\R$
determines a cocycle $f^{(\,\cdot\,)}(\,\cdot\,):\Z\times X\to \R$
given by
\[f^{(m)}(x)=\left\{\begin{array}{ccc}
f(x)+f(Tx)+\ldots+f(T^{m-1}x) & \mbox{if} & m>0\\ 0 & \mbox{if} &
m=0\\ -\left(f(T^mx)+\ldots+f(T^{-1}x)\right)  & \mbox{if} & m<0.
\end{array}\right.\]

Denote by $\lambda$ Lebesgue measure on $\R$. If $f:X\to\R$ is a
strictly positive $L^1$--function, then by $T^f=(T^f_t)_{t\in\R}$
we will mean the corresponding special flow under $f$ (see e.g.\
\cite{Co-Fo-Si}, Chapter 11) acting on $(X^f,{\cal B}^f,\mu^f)$,
where $X^f=\{(x,s)\in X\times \R:\:0\leq s<f(x)\}$ and ${\cal
B}^f$ $(\mu^f)$ is the restriction of ${\cal B}\otimes{\cal
B}(\R)$ $(\mu\otimes\lambda)$ to $X^f$. Under the action of the
flow $T^f$ each point in $X^f$ moves vertically at unit speed, and
we identify the point $(x,f(x))$ with $(Tx,0)$. More precisely, if
$(x,s)\in X^f$ then
\[T^f_t(x,s)=(T^nx,s+t-f^{(n)}(x)),\]
where $n\in\Z$ is a unique number such that
\[f^{(n)}(x)\leq s+t<f^{(n+1)}(x).\]

We denote by $\T$ the circle group $\R/\Z$ which we will
constantly identify with the interval $[0,1)$ with addition mod
$1$.  For a real number $t$ denote by $\{t\}$ its fractional part
and by $\|t\|$ its distance to the nearest integer number.  For an
irrational $\alpha\in\T$ denote by $(q_n)$ its sequence of
denominators (see  e.g.\ \cite{Ch}), that is we have
\begin{equation}\label{ulla}
\frac{1}{2q_nq_{n+1}}<\left|\alpha-\frac{p_n}{q_n}\right|<\frac{1}{q_nq_{n+1}},
\end{equation}
where
\[\begin{array}{ccc}
q_0=1, & q_1=a_1, & q_{n+1}=a_{n+1}q_n+q_{n-1}\\
p_0=0, & p_1=1,  & p_{n+1}=a_{n+1}p_n+p_{n-1}
\end{array}\]
and $[0;a_1,a_2,\dots]$ stands for the continued fraction
expansion of $\alpha$.  We say that $\alpha$ has {\em bounded
partial quotients} if the sequence $(a_n)$ is bounded. If
$C=\sup\{a_n:n\in\N\}+1$ then
\[\frac{1}{2Cq_n}<\frac{1}{2q_{n+1}}<\|q_n\alpha\|<\frac{1}{q_{n+1}}<\frac{1}{q_n}\]
for each $n\in\N$.

\section{Construction}\label{constr}
In this section, using the procedure of gluing of flows which was
described by Blokhin in \cite{Bl}, we will construct the flow
$(\varphi_t)_{t\in\R}$ that was announced in Introduction.

Let $\alpha\in\R$ be an irrational number. We denote by
$(\psi_t)_{t\in\R}$ the linear flow
$\psi_t(x_1,x_2)=(x_1+t\alpha,x_2+t)$ of the torus $\T^2$.

\begin{figure}[h]
\begin{center}
\vspace{5.5cm} \includegraphics{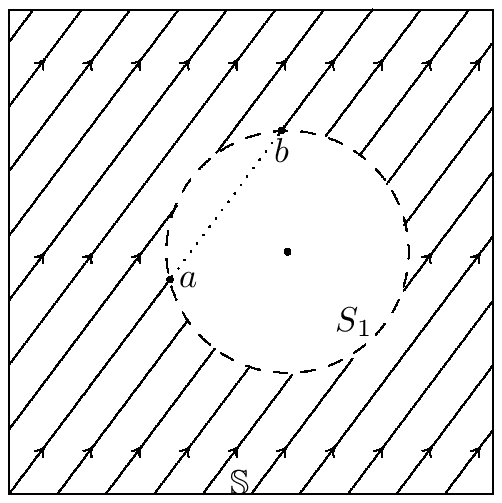}
\end{center}\caption{}
\end{figure}

Let us cut out (from the torus) a  disk $D$ which is disjoint from
the circle $\s=\{(x,0)\in\T^2:x\in[0,1)\}$ and intersects the
segment $O_{[0,1]}=\{(\alpha t,t):t\in[0,1]\}$. We will denote by
$S_1$ the circle which bounds $D$. Let $a$ and $b$ be points of
$S_1$ that lie on the segment $O_{[0,1]}$ (see Fig.1).

\begin{figure}[h]
\begin{center}
\vspace{5.5cm} \includegraphics{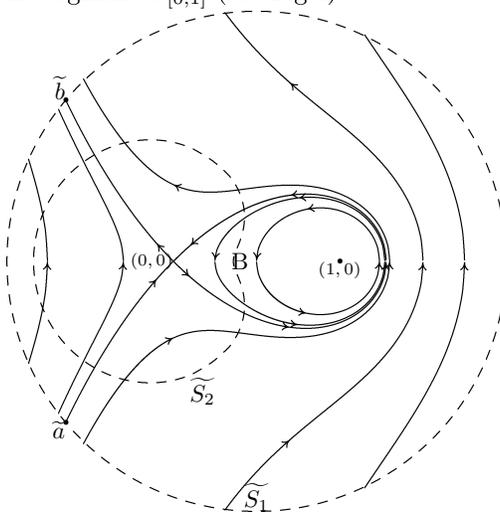}
\end{center}\caption{The phase
portrait for the Hamiltonian system $H(x,y)=\frac{1}{2}
e^{2x}(y^2+(x-1)^2)$; the portrait is the same as
for~(\ref{rowroz})}
\end{figure}

Now let us consider the flow $(\widetilde{\psi}_t)_{t\in\R}$ on
the disk
$\widetilde{D}=\{(x,y)\in\R^2:\:(x-1/2)^2+y^2\leq(3/2)^2\}$ given
by the system of equations (on $\R^2\setminus\{(0,0)\}$)
\begin{equation}\label{rowroz}
\left\{\begin{array}{rcl}
\frac{dx}{dt}&=&\frac{-y}{x^2+y^2}\\
\frac{dy}{dt}&=&\frac{x(x-1)+y^2}{x^2+y^2}.
\end{array}\right.
\end{equation}
By the Liouville theorem, $(\widetilde{\psi}_t)_{t\in\R}$
preserves the measure $e^{2x}(x^2+y^2)\,dx\,dy$. Let
$\widetilde{S_1}$ be the boundary of $\widetilde{D}$. The flow
$(\widetilde{\psi}_t)$ has a singularity at $(0,0)$ and a fixed
point $(1,0)$ which is a center (see Fig.\ 2). Moreover the set
\[B=\{(x,y)\in\widetilde{D}:e^{2x}(y^2+(x-1)^2)<1,\;x>0\}.\]
consists of periodic orbits. Let $\widetilde{a}$ and
$\widetilde{b}$ be the points of intersection of $\widetilde{S_1}$
and the separatrices of $(0,0)$. By Lemma~1 in \cite{Bl}, there
exists a $C^{\infty}$--diffeomorphism $g:S_1\to \widetilde{S_1}$
such that $g(a)=\widetilde{a}$, $g(b)=\widetilde{b}$, and there
exist a $C^{\infty}$--flow $(\varphi_t)_{t\in\R}$ on
$M=(\T\setminus D)\cup_g \widetilde{D}$ and a
$C^{\infty}$--measure $\mu$ on $(\T\setminus D)\cup_g
\widetilde{D}$ such that
\begin{itemize}
\item $(\varphi_t)_{t\in\R}$ preserves $\mu$,
\item the flow $(\varphi_t)_{t\in\R}$ restricted to $\T\setminus
D$ is equal to $(\psi_t)_{t\in\R}$,
\item the flow $(\varphi_t)_{t\in\R}$ restricted to $\widetilde{D}$ is equal to
$(\widetilde{\psi}_t)_{t\in\R}$.
\end{itemize}
$M$ splits into two  $(\varphi_t)_{t\in\R}$--invariant sets $B$
and $A=M\setminus B$ such that $B$ consists of periodic orbits and
$A$ is an ergodic component of positive measure. Moreover, the
flow $(\varphi_t)_{t\in\R}$ on $A$ can be represented as the
special flow built over the rotation $Tx=x+\alpha$ and under a
function $f:\T\to\R$ which is of class $C^{\infty}$ on
$\T\setminus\{0\}$. Of course, $f(x)$ is the first return time to
$\s$ of the point $x\in\s\cong\T$. We will prove that
$f:(0,1)\to\R$ can be extended to a $C^{\infty}$ function on
$[0,1]$, i.e.\ $D^nf$ possesses limits at $0$ and $1$ for any
$n\geq 0$. Moreover, we will show that
$\lim_{x\to0^+}f(x)>\lim_{x\to1^-}f(x)$. To prove it we will need
an auxiliary simple lemma.

\begin{Lemma} Let $U\subset\C$ be an open disk with center at $0$ and $h:U\to\C$
be an analytic function such that $h(z)\neq 0$ for $z\in U$. Let
us consider the differential equation
\[\frac{dz}{dt}=\frac{i}{zh(z)}\]
on $U\setminus\{0\}$. Then there exists an open disk
$\widetilde{U}\subset U$ containing $0$ and  a biholomorphic map
$\xi:\widetilde{U}\to\xi(\widetilde{U})$ such that $\xi(0)=0$ and
\[\frac{d\omega}{dt}=1/\omega\]
on $\xi(\widetilde{U})\setminus\{0\}$, where
$\omega=\sqrt{2}\xi(z)$.
\end{Lemma}
\begin{proof}
Let $H:U\to\C$ be an analytic function such that $H'(z)=-izh(z)$
and $H(0)=0$. Since $H'(0)=0$ and $H''(0)=-ih(0)\neq 0$, there
exists an open disk $\widetilde{U}\subset U$ containing $0$ and  a
biholomorphic map $\xi:\widetilde{U}\to\xi(\widetilde{U})$ such
that $H(z)=(\xi(z))^2$ for $z\in\widetilde{U}$. Put
$\omega=\sqrt{2}\xi(z)$, $z\in\tilde{U}$. Then
$\omega^2/2=(\xi(z))^2=H(z)$, and consequently
\[\frac{d\omega}{dt}\omega=H'(z)\frac{dz}{dt}=\frac{iH'(z)}{zh(z)}=1.\]
\end{proof}

Of course, the equation (\ref{rowroz}) can be written as
$\frac{dz}{dt}=i(z-1)/z$. By Lemma~1 (with $h(z)=\frac1{z-1}$),
there exist an open disk $0\in V\subset \C$ and a biholomorphic
map $F:V\to F(V)$ such that the flow $(F^{-1}\circ\varphi_t\circ
F)_{t\in\R}$ on $V$ is determined by the equation
$\frac{d\omega}{dt}=1/\omega$, i.e.\ by
\begin{equation}\label{rowrozhip}
\left\{\begin{array}{rcl}
\frac{dx}{dt}&=&\frac{x}{x^2+y^2}\\
\frac{dy}{dt}&=&\frac{-y}{x^2+y^2}.
\end{array}\right.
\end{equation}

\begin{figure}[h]
\begin{center}
\vspace{5.5cm} \includegraphics{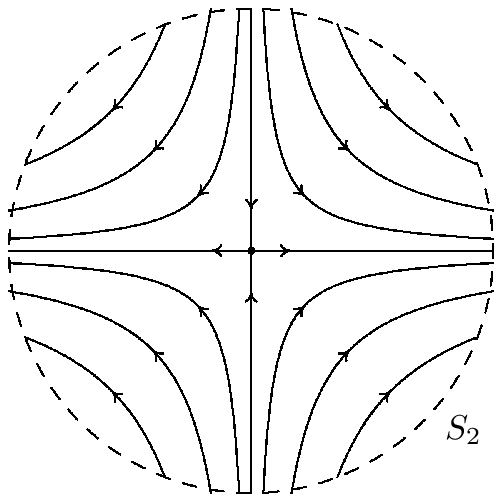}
\end{center}\caption{}
\end{figure}

 The trajectories of this flow are presented on Fig.\ 3. Denote
by $S_2=\{re^{it}:t\in[0,2\pi]\}$ a circle which is contained in
$V$. Let $\tau:S_2\to\R$ be the function of first return time
(counted forward or backward and staying inside $S_2$) to $S_2$.
It is easy to check that
\[\tau(re^{it})=-r^2\cos(2t)\]
which is of class $C^{\infty}$ (indeed, $\frac{d}{dt}(\omega^2)=2$
and the first return time satisfies $|\omega^2(t)|^2=r^4$). Let
$\widetilde{S_2}:=F(S_2)$ and
$\widetilde{\tau}:\widetilde{S_2}\to\R$  be the function of the
first return time (counted forward or backward inside $F(V)$) to
$\widetilde{S_2}$. Then $\widetilde{\tau}$ is also of class
$C^{\infty}$. Consequently, $f:(0,1)\to\R$ can be extended to a
$C^{\infty}$ function on $[0,1]$ and
\[\lim_{x\to0^+}f(x)-\lim_{x\to1^-}f(x)\geq \tau_0,\]
where $\tau_0$ is the time of the first positive return of the
point $0$ to itself via the separatrice which starts and stops at
$0$. In this way we constructed a $C^{\infty}$--flow with one
singular point on the torus which
\begin{itemize}
\item preserves a $C^{\infty}$--measure,
\item possesses two invariant subsets $A$ and $B$: $A$ is an
ergodic component of positive measure and $B$ consists of periodic
orbits,
\item the flow on $A$ is measure--theoretically isomorphic to a
special flow $T^f$, where $T$ is the rotation by $\alpha$ and
$f:\T\to\R$ is $C^{\infty}$ function on $\T\setminus\{0\}$ and
$\lim_{x\to0^+}f(x)\neq\lim_{x\to0^-}f(x)$.
\end{itemize}

By cutting out the disk with the center at $(1,0)$ and of radius
$1/2$ (it intersects $A$ but does not contain the point $(0,0)$)
from the flow $(\varphi_t)_{t\in\R}$ and gluing a M\"{o}bius strip
endowed with the flow considered by Blokhin in \cite[\S3]{Bl} we
can obtain a $C^\infty$ flow $(u_t)_{t\in\R}$ with one singularity
on a non-orientable surface of Euler characteristic -1 such that
$(u_t)$ is isomorphic to the action of $(\varphi_t)$  on the
component $A$.

\section{Joinings}
Assume that $\cs=\st$ is a flow on a standard probability space
$\xbm$. By that we mean always a so called {\em measurable flow},
i.e.\ we require in particular that the map $\R\ni t\to\langle
f\circ S_t,g\rangle\in\C$ is continuous for each $f,g\in L^2\xbm$.
Assume moreover that $\cs$ is ergodic and let $\ct=\ttt$ be
another ergodic flow defined on $\ycn$. By a {\em joining} between
$\cs$ and $\ct$ we mean any probability $(S_t\times
T_t)_{t\in\R}$--invariant measure on $(X\times Y,{\cal B}\otimes
{\cal C})$ whose projections on $X$ and $Y$ are equal to $\mu$ and
$\nu$ respectively. The set of joinings between $\cs$ and $\ct$ is
denoted by $J(\cs,\ct)$. The subset of ergodic joinings is denoted
by $J^e(\cs,\ct)$. Ergodic joinings are exactly extremal points in
the simplex $J(\cs,\ct)$. Let $\{A_n:n\in\N\}$ and
$\{B_n:n\in\N\}$ be two countable families in $\mathcal{B}$ and
$\mathcal{C}$ respectively which are dense in $\mathcal{B}$ and
$\mathcal{C}$ for the (pseudo--)metrics
$d_{\mu}(A,B)=\mu(A\triangle B)$ and $d_{\nu}(A,B)=\nu(A\triangle
B)$ respectively. Let us consider the metric $d$ on $J(\cs,\ct)$
defined by
\[d(\rho,\rho')=\sum_{m,n\in\N}\frac{1}{2^{m+n}}|\rho(A_n\times B_m)-\rho'(A_n\times B_m)|.\]
Endowed with corresponding to $d$ topology, to which we will refer
as the weak topology, the set $J(\cs,\ct)$ is compact.

Suppose that $\mathcal{A}\subset\mathcal{B}$ is a {\em factor} of
$\cs$, i.e.\  $\mathcal{A}$ is an $\cs$--invariant
sub--$\sigma$--algebra. Denote by $\mu\otimes_{\mathcal{A}}\mu\in
J(\cs,\cs)$ the {\em relatively independent joining} of the
measure $\mu$ over the factor $\mathcal{A}$, i.e.\
$\mu\otimes_{\mathcal{A}}\mu\in J(\cs,\cs)$ is defined by
\[(\mu\otimes_{\mathcal{A}}\mu)(D)=
\int_{X/\mathcal{A}}(\mu_{\overline{x}}\otimes\mu_{\overline{x}})(D)\,d\overline{\mu}(\overline{x})\]
for $D\in\mathcal{B}\otimes\mathcal{C}$, where
$\{\mu_{\overline{x}}:\overline{x}\in X/\mathcal{A}\}$ is the
disintegration of the measure $\mu$ over the factor $\mathcal{A}$
and $\overline{\mu}$ is the image of $\mu$ by the factor map $X\to
X/\mathcal{A}$.

For every $t\in\R$ by $\mu_{S_t}\in J(\cs,\cs)$ we will denote the
graph joining determined by $\mu_{S_t}(A\times B)=\mu(A\cap
S_{-t}B)$ for $A,B\in\mathcal{B}$. Then $\mu_{S_t}$ is
concentrated on the graph of $S_t$ and $\mu_{S_t}\in
J^e(\cs,\cs)$.

Let $(t_n)_{n\in\N}$ be a sequence of real numbers such that
$t_n\to+\infty$. We say that a flow $\cs$ on $\xbm$ is {\em rigid}
along $(t_n)$ if
\begin{equation}\label{defszt}
\mu(A\cap S_{-t_n}A)\to \mu(A)
\end{equation}
for every $A\in\mathcal{B}$, or, equivalently, $\mu_{S_{t_n}}\to
\mu_{Id}$ weakly in $J(\cs,\cs)$. In particular, a factor
$\mathcal{A}\subset\mathcal{B}$ of $\cs$ is rigid along $(t_n)$ if
the convergence (\ref{defszt}) holds for every $A\in\mathcal{A}$.
 It is well known that a flow is mildly mixing iff it has no non--trivial
rigid factor (see~\cite{Fu-We,Sch-Wa}).

\begin{definition}
A flow $\cs$ on $\xbm$ is called {\em partially rigid} along
$(t_n)$ if there exists $0<u\leq 1$ such that
\[\liminf_{n\to\infty}\mu(A\cap S_{-t_n}A)\geq u\mu(A)\;\;\text{ for every }\;\;A\in\mathcal{B},\]
or, equivalently, every weak limit point $\rho$ of the sequence
$(\mu_{S_{t_n}})_{n\in\N}$ in $J(\cs,\cs)$ satisfies
$\rho(\Delta)\geq u$, where $\Delta=\{(x,x)\in X\times X:x\in
X\}$.
\end{definition}

The proof of the following proposition is the same as in the case
of measure--preserving transformations and can be found in
\cite{Le-Pa}.

\begin{Prop}\label{prlepa}
Let $\cs$ be an ergodic flow on $\xbm$. Suppose that
$\mathcal{A}\subset\mathcal{B}$ is a non--trivial rigid factor of
$\cs$. Then there exist a factor $\mathcal{A}'\supset\mathcal{A}$
of $\cs$ and a rigidity sequence $(t_n)$ for $\mathcal{A}'$ such
that $\mu_{S_{t_n}}\to \mu\otimes_{\mathcal{A}'}\mu$ weakly in
$J(\cs,\cs)$.
\end{Prop}

Recall that in general the notions of (absence of) partial
rigidity and mild mixing are not related. For example, the Chacon
transformation is partially rigid (see e.g.\ \cite{Pr-Ry}) and
mildly mixing. On the other hand the Cartesian product of a mixing
transformation and a rigid transformation is not mildly mixing and
has no partial rigidity. Under some additional strong assumption
we have however the following.

\begin{Lemma}\label{jmm}
Let $\cs$ be an ergodic flow on $\xbm$ which is a finite extension
of each of its non--trivial factors.  Then if the flow $\cs$ is
not partially rigid then it is mildly mixing.
\end{Lemma}

\begin{proof}
Suppose, contrary to our claim, that there exists a non--trivial
factor $\mathcal{A}$ of $\cs$ which is rigid. By
Proposition~\ref{prlepa} there exist a factor
$\mathcal{A}'\supset\mathcal{A}$ and a rigidity sequence $(t_n)$
for $\mathcal{A}'$ such that
\begin{equation}\label{zbdo}
\mu_{S_{t_n}}\to \mu\otimes_{\mathcal{A}'}\mu\text{ weakly in
}J(\cs,\cs).
\end{equation}
Since $\cs$ is ergodic and it  is a finite extension of
$\cs|_{\mathcal{A}'}$, there exists a natural number $k$ such that
every fiber measure $\mu_{\overline{x}}$, $\overline{x}\in
X/\mathcal{A}'$ is atomic with $k$ atoms each of measure $1/k$,
where $\{\mu_{\overline{x}}:\overline{x}\in X/\mathcal{A}'\}$ is
the disintegration of the measure $\mu$ over the factor
$\mathcal{A}'$. Then
\[(\mu\otimes_{\mathcal{A}'}\mu)(\Delta)=
\int_{X/\mathcal{A}}(\mu_{\overline{x}}\otimes\mu_{\overline{x}})(\Delta)\,d\overline{\mu}(\overline{x})=
\int_{X/\mathcal{A}}\frac{1}{k}\,d\overline{\mu}(\overline{x})=\frac{1}{k},\]
which, in view of (\ref{zbdo}), gives the partial rigidity of
$\cs$ and we obtain a contradiction.
\end{proof}

\section{Ratner property}\label{ratprop}
In this section we introduce and study  a condition which emulates
the Ratner condition from \cite{Rat}.
\begin{definition} \label{defrat}(cf.\ \cite{Rat,Th})
Let $(X,d)$ be a $\sigma$--compact metric space, $\mathcal{B}$ be
the $\sigma$--algebra of Borel subsets of $X$, $\mu$ a Borel
probability measure on $(X,d)$ and let $(S_t)_{t\in\R}$ be a  flow
on the space $(X,{\cal B},\mu)$. Let $P\subset\R\setminus\{0\}$ be
a finite subset and $t_0\in\R\setminus\{0\}$. The flow
$(S_t)_{t\in\R}$ is said to have {\em the property $\rat(t_0,P)$}
if for every $\vep>0$ and $N\in\N$ there exist
$\kappa=\kappa(\vep)>0$, $\delta=\delta(\vep,N)>0$ and a subset
$Z=Z(\vep,N)\in\mathcal{B}$ with $\mu(Z)>1-\vep$  such that if
$x,x'\in Z$, $x'$ is not in the orbit $x$ and $d(x,x')<\delta$,
then there are $M=M(x,x')$, $L=L(x,x')\geq N$ such that
$L/M\geq\kappa$ and there exists $p=p(x,x')\in P$ such that
\[\frac{\# \{n\in\Z\cap[M,M+L]:d(S_{nt_0}(x),S_{nt_0+p}(x'))<\vep\}}{L}>1-\vep.\]
Moreover, we say that $(S_t)_{t\in\R}$ has {\em the property
$\rat(P)$} if the set of all $s\in\R$ such that the flow
$(S_t)_{t\in\R}$ has the $\rat(s,P)$--property is uncountable.
\end{definition}

\begin{remark}
In the original definition of M.\ Ratner, for $t_0\neq 0$,
$P=\{-t_0,t_0\}$. In our situation a priori there is no relation
between $t_0$ and $P$. Analysis similar to that in the proof of
Theorem~2 in \cite{Rat} shows that $\rat(t_0,P)$ and $\rat(P)$
properties are invariant under measure--theoretic isomorphism.
\end{remark}

We now prove an extension of Theorem~3 in \cite{Rat} that brings
important information about ergodic joinings with flows satisfying
the $\rat(P)$--property.

\begin{Th}\label{prpr}
Let $(X,d)$ be a $\sigma$--compact metric space, $\mathcal{B}$ be
the $\sigma$--algebra of Borel subsets of $X$ and $\mu$ a
probability measure on $(X,\mathcal{B})$. Let
$P\subset\R\setminus\{0\}$ be a nonempty finite set. Assume that
$(S_t)_{t\in\R}$ is an ergodic flow on $(X,{\cal B},\mu)$ such
that every automorphism $S_p:(X,{\cal B},\mu)\to(X,{\cal B},\mu)$
for $p\in P$ is ergodic. Suppose that $(S_t)_{t\in\R}$ satisfies
the $\rat(P)$--property. Let $(T_t)_{t\in\R}$ be an ergodic flow
on $(Y,{\cal C},\nu)$ and let $\rho$ be an ergodic joining of
$(S_t)_{t\in\R}$ and $(T_t)_{t\in\R}$. Then either
$\rho=\mu\otimes\nu$, or $\rho$ is a finite extension of $\nu$.
\end{Th}

\begin{remark}\label{finext}
Let $\cs$ be an ergodic flow on $\xbm$. Assume that for each
ergodic flow $\ct$ acting on $\ycn$ an arbitrary ergodic joining
$\rho$ of $\cs$  with  $\ct$ is either the product measure or
$\rho$ is a finite extension of $\mu$. Then $\cs$ is a finite
extension of each of its non--trivial factors. Indeed, suppose
that $\mathcal{A}\subset\mathcal{B}$ is a non--trivial factor. Let
us consider the factor action $\cs|_{\mathcal{A}}$ on
$(X/\mathcal{A},\mathcal{A},\overline{\mu})$ and the natural
joining $\mu_{\mathcal{A}}\in J(\cs,\cs|_{\mathcal{A}})$
determined by $\mu_{\mathcal{A}}(B\times A)=\mu(B\cap A)$ for all
$B\in\mathcal{B}$ and $A\in\mathcal{A}$. Clearly, the action
$\cs\times(\cs|_{\mathcal{A}})$ on $(X\times
(X/\mathcal{A}),\mathcal{B}\otimes\mathcal{A},\mu_{\mathcal{A}})$
is isomorphic (via the projection on $\xbm$) to the action of
$\cs$. Since the measure $\mu_{\mathcal{A}}$  is not the product
measure, by assumptions, the action
$\cs\times(\cs|_{\mathcal{A}})$ on $(X\times
(X/\mathcal{A}),\mathcal{B}\otimes\mathcal{A},\mu_{\mathcal{A}})$
is a finite extension of $\cs|_{\mathcal{A}}$.
\end{remark}

To prove Theorem~\ref{prpr} we will need two ingredients. The
proof of the following lemma is contained in the proof of
Theorem~3 in \cite{Rat}.

\begin{Lemma}\label{lemrat}
Let $(S_t)_{t\in\R}$ and $(T_t)_{t\in\R}$ be  ergodic flows acting
on $\xbm$ and $\ycn$  respectively and let $\rho\in J^e(\cs,\ct)$.
Suppose that there exists $U\in\mathcal{B}\otimes\mathcal{C}$ with
$\rho(U)>0$ and $\delta>0$ such that   if $(x,y)\in U$, $(x',y)\in
U$ then either $x$ and $x'$ are in the same orbit or
$d(x,x')\geq\delta$. Then $\rho$ is a finite extension of $\nu$.
\end{Lemma}

The following simple fact will be used in the proof of
Theorem~\ref{prpr} and in the remainder of the paper.

\begin{remark}\label{uwaerg}
Notice that if
\[\left|\frac{1}{M}\sum_{n=0}^{M-1}\chi_A(T^nx)-\mu(A)\right|<\vep\text{
and
}\left|\frac{1}{M+L+1}\sum_{n=0}^{M+L}\chi_A(T^nx)-\mu(A)\right|<\vep\]
then
\[\left|\frac{1}{ L}\sum_{n=M}^{M+L}\chi_A(T^nx)-\mu(A)\right|<2\vep\left(1+\frac{M}{L}\right).\]
\end{remark}

The proof of Theorem~\ref{prpr}, presented below, is much the same
as the proof of Theorem~10 in \cite{Th}.

\begin{proofpr}
Suppose that  $\rho\in J^e(\cs,\ct)$ and $\rho\neq\mu\otimes\nu$.
Since the flow $(S_t\times T_t)_{t\in\R}$ is ergodic on $(X\times
Y,\rho)$, we can find $t_0\neq 0$ such that the automorphism
$S_{t_0}\times T_{t_0}:(X\times Y,\rho)\to(X\times Y,\rho)$ is
ergodic and the flow $(S_{t})_{t\in\R}$  has the
$\rat(t_0,P)$--property. To simplify notation we assume that
$t_0=1$.

 Since the ergodicity of $S_p$
implies disjointness of $S_{p}$  from the identity, for every
$p\in P$ there exist closed subsets $A_p\subset X$, $B_p\subset Y$
such that
\[\rho(S_{-p}A_p\times B_p)\neq \rho(A_p\times B_p).\]
 Let
\begin{equation}\label{mini}
0<\vep:=\min\{|\rho(S_{-p}A_p\times B_p)-\rho(A_p\times B_p)|:p\in
P\}.
\end{equation}
Next choose $0<\vep_1<\vep/8$ such that $\mu(A_p^{\vep_1}\setminus
A_p)<\vep/2$ for $p\in P$, where $A^{\vep_1}=\{z\in
X:d(z,A)<\vep_1\}$. We have
\begin{equation}\label{szaep1}
\begin{aligned} |\rho(A_p\times
B_p)-\rho(A_p^{\vep_1}\times B_p)| &=&  \rho(A_p^{\vep_1}\times
B_p\setminus A_p\times
B_p)\;\;\;\;\;\;\;\;\;\;\;\;\;\;\;\;\;\;\;\;\;\;\;\;\;\;\\& \leq&
\rho((A_p^{\vep_1}\setminus A_p)\times Y)
=\mu(A_p^{\vep_1}\setminus A_p)<\vep/2
\end{aligned}
\end{equation}
and similarly
\begin{equation*}\label{szaep2}
|\rho(S_{-p}A_p\times B_p)-\rho(S_{-p}(A_p^{\vep_1})\times
B_p)|<\vep/2
\end{equation*}
 for any $p\in P$.

Let $\kappa:=\kappa(\vep_1)(>0)$. By the ergodic theorem together
with Remark~\ref{uwaerg}, there exist a measurable set $U\subset
X\times Y$ with $\rho(U)>3/4$ and $N\in\N$ such that if $(x,y)\in
U$, $p\in P$, $m\geq N$ and $l/m\geq \kappa$ then
\begin{equation}\label{some2}
\left|\frac{1}{l}\sum_{k=m}^{m+l}\chi_{A_p^{\vep_1}\times
B_p}(S_{k}x,T_{k}y)-\rho(A_p^{\vep_1}\times
B_p)\right|<\frac{\vep}{8},
\end{equation}
\begin{equation}\label{some3}
\left|\frac{1}{l}\sum_{k=m}^{m+l}\chi_{S_{-p}A_p\times
B_p}(S_{k}x,T_{k}y)-\rho(S_{-p}A_p\times
B_p)\right|<\frac{\vep}{8}
\end{equation}
and similar inequalities hold for $A_p\times B_p$ and
$S_{-p}(A_p^{\vep_1})\times B_p$.

Next, by the property $\rat(1,P)$, we obtain relevant
$\delta=\delta(\vep_1,N)>0$ and $Z=Z(\vep_1,N)\in\mathcal{B}$,
$\mu(Z)>1-\vep_1$.

Now assume that $(x,y)\in U$, $(x',y)\in U$, $x,x'\in Z$ and $x'$
is not in the orbit of $x$. We claim that $d(x,x')\geq \delta$.
Suppose that, on the contrary, $d(x,x')<\delta$. Then, by the
property $\rat(1,P)$, there exist $M=M(x,x')$, $L=L(x,x')\geq N$
with $L/M\geq\kappa$ and $p=p(x,x')\in P$ such that
$(\#K_p)/L>1-\vep_1$, where
\[K_p=\{n\in\Z\cap[M,M+L]:d(S_{n}(x),S_{n+p}(x'))<\vep_1\}.\]
If $k\in K_p$ and $S_{k+p}x'\in A_p$, then $S_{k}x\in
A_p^{\vep_1}$. Hence
\begin{equation}\label{kicha1}
 \begin{aligned}
\lefteqn{\frac{1}{L}\sum_{k=M}^{M+L}\chi_{S_{-p}A_p\times
B_p}(S_{k}x',T_{k}y)}\\
&\leq \frac{\#(\Z\cap[M,M+L]\setminus K_p)}{L}+\frac{1}{L}\sum_{k\in K_p}\chi_{A_p\times B_p}(S_{k+p}x',T_{k}y)\\
& \leq \vep/8+\frac{1}{L}\sum_{k=M}^{M+L}\chi_{A_p^{\vep_1}\times
B_p}(S_{k}x,T_{k}y).
\end{aligned}
\end{equation}
Now from (\ref{some3}), (\ref{kicha1}), (\ref{some2}) and
(\ref{szaep1}) it follows that
\begin{eqnarray*}
\rho(S_{-p}A_p\times B_p) &\leq &
\frac{1}{L}\sum_{k=M}^{M+L}\chi_{S_{-p}A_p\times
B_p}(S_{k}x',T_{k}y)+\vep/8\\
& \leq &
\vep/4+\frac{1}{L}\sum_{k=M}^{M+L}\chi_{A_p^{\vep_1}\times
B_p}(S_{k}x,T_{k}y)\\
& < & \vep/2+\rho(A_p^{\vep_1}\times B_p)\leq \vep+\rho(A_p\times
B_p).
\end{eqnarray*}
Applying similar arguments we get
\[\rho(A_p\times B_p)<
\vep+\rho(S_{-p}A_p\times B_p).\] Consequently,
\[|\rho(A_p\times B_p)-\rho(S_{-p}A_p\times B_p)|<\vep,\]
contrary to (\ref{mini}).

In summary, we found a measurable set $U_1=U\cap(Z(\vep_1,N)\times
Y)$ and $\delta(\vep_1,N)>0$ such that $\rho(U_1)>3/4-\vep_1>1/2$
 and if $(x,y)\in U_1$,  $(x',y)\in U_1$ then either $x$ and
$x'$ are in the same orbit or $d(x,x')\geq\delta(\vep_1,N)$. Now
an application of Lemma~\ref{lemrat} completes the proof.\hfill
$\Box$
\end{proofpr}
We end up this section with  a general lemma that gives a
criterion that allows one to prove the $\rat(P)$--property for
special flows built over irrational rotations on the circle and
under bounded and bounded away from zero measurable functions.

While dealing with special flows over irrational rotations on
$\T^f$ we will always consider the induced metric from the metric
defined on $\T\times\R$ by $d((x,s),(y,t))=\|x-y\|+|s-t|$.

\begin{Lemma}\label{lemfundsp}
Let $P\subset\R\setminus\{0\}$ be a nonempty finite subset. Let
$T:\T\to\T$ be an ergodic rotation and let $f:\T\to\R$ be a
bounded positive measurable function which is bounded away from
zero. Assume that for every $\vep>0$ and $N\in\N$ there exist
$\kappa=\kappa(\vep)>0$ and $0<\delta=\delta(\vep,N)<\vep$ such
that if $x,y\in\T$, $0<\|x-y\|<\delta$, then there are natural
numbers $M=M(x,y)\geq N$, $L=L(x,y)\geq N$ such that
$L/M\geq\kappa$ and there exists $p=p(x,y)\in P$ such that
\[\frac{1}{L+1}\#\left\{M\leq n\leq M+L:|f^{(n)}(x)-f^{(n)}(y)+p|<\vep\right\}>1-\vep.\]
Suppose that $\gamma\in\R$ is a positive number such that the
instance automorphism $T^f_\gamma:\T^f\to\T^f$ is ergodic. Then
the special flow $T^f$ has the $\rat(\gamma,P)$--property.
\end{Lemma}

\begin{proof}
Let $c$, $C$ be positive numbers such that $0<c\leq f(x)\leq C$
for every $x\in\T$. Let $\mu$ stand for Lebesgue measure on $\T$.

Fix $0<\vep$ and $N\in\N$.  Put
\[\vep_1=\min\left(\frac{c\vep}{8(\gamma+C)},\frac{\vep}{16}\right).\]
Take $\kappa'=\kappa(\vep_1)$ and let
$\kappa:=\frac{c}{C}\kappa'$. Let us consider the set
\[X(\vep):=\left\{(x,s)\in\T^f:\frac{\vep}{8}<s<f(x)-\frac{\vep}{8}\right\}.\]
Since $\mu^f(X(\vep)^c)=\vep/4$ and $T^f_\gamma$ is ergodic, there
exists $N(\vep)\in\N$ such that $\mu^f(Z(\vep)^c)<\vep$, where
$Z(\vep)$ is the set of all $(x,s)\in\T^f$ such that
\begin{equation}\label{wsgeba}
\left|\frac{1}{n}\#\{0\leq k<n:T^f_{k\gamma}(x,s)\notin
X(\vep)\}-\frac{\vep}{4}\right|<\frac{\kappa}{1+\kappa}\frac{\vep}{8}
\end{equation}
for each $n\geq N(\vep)$. Take
$\delta=\delta(\vep_1,2\gamma\max(N(\vep),N)/c)<\vep_1$.
 Let us consider a pair of
points $(x,s),(y,s')\in Z(\vep)$ such that
$0<d((x,s),(y,s'))<\delta$ and $x\neq y$. By assumption, there are
natural numbers $M'=M(x,y), L'=L(x,y)\geq
2\gamma\max(N(\vep),N)/c$ such that $L'/M'\geq\kappa'$ and there
exists $p=p(x,y)\in P$ such that \[\frac{\#A'}{L'+1}>1-\vep_1,\]
where $A'=\left\{M'\leq n\leq
M'+L':|f^{(n)}(x)-f^{(n)}(y)+p|<\vep_1\right\}$. Then
\begin{equation}\label{nierapr}
\frac{\#A''}{L'}>1-4\vep_1, \;\text{ where }\;A''=\left\{M'\leq
n<M'+L':n\in A',n+1\in A'\right\}.
\end{equation}
 Put
\[M:=\frac{f^{(M')}(x)-s}{\gamma}\;\;\text{ and }\;\;L:=\frac{f^{(L')}(T^{M'}x)}{\gamma}.\]
Then
\[\frac{L}{M}=\frac{f^{(L')}(T^{M'}x)}{f^{(M')}(x)-s}\geq
\frac{c}{C}\frac{L'}{M'}\geq \kappa.\] But $s\leq f(x)$ and
 $M',L'\geq 2\gamma\max(N(\vep),N)/c$, so
\[M=\frac{f^{(M')}(x)-s}{\gamma}\geq \frac{f^{(M'-1)}(Tx)}{\gamma}\geq\frac{c(M'-1)}{\gamma}\geq
 \frac{cM'}{2\gamma}\geq\max(N(\vep),N)\] and
\begin{equation}\label{nierl}
L=\frac{f^{(L')}(T^{M'}x)}{L'}\frac{L'}{\gamma}\geq
c\frac{L'}{\gamma}>N.
\end{equation}
Now $M\geq N(\vep)$, $L/M\geq\kappa$, $(x,s)\in Z(\vep)$ (that is
$(x,s)$ satisfies (\ref{wsgeba})) so, by Remark~\ref{uwaerg}, we
have
\begin{equation}\label{geba}
\frac{1}{L}\#\{M\leq k<M+L:T^f_{k\gamma}(x,s)\notin
X(\vep)\}<\frac{\vep}{2}.
\end{equation}
Suppose that $M\leq k<M+L$. Then
$k\gamma+s\in[f^{(M')}(x),f^{(M'+L')}(x))$ and there exists a
unique $M'\leq m_k<M'+L'$ such that
$k\gamma+s\in[f^{(m_k)}(x),f^{(m_k+1)}(x))$. Suppose that
\[k\in B:=\{M\leq j<M+L:T^f_{j\gamma}(x,s)\in X(\vep)\text{ and }m_j\in
A''\}.\] Then
\[f^{(m_k)}(x)+\vep/8< s+k\gamma<f^{(m_k+1)}(x)-\vep/8.\]
Since $m_k\in A''$ and $|s'-s|<\delta<\vep_1$, we have
\begin{eqnarray*}
s'+k\gamma+p & = & (s+k\gamma)+(s'-s)+p <
 f^{(m_k+1)}(x)+p-\vep/8+\delta\\
& =& f^{(m_k+1)}(y)+(f^{(m_k+1)}(x)-f^{(m_k+1)}(y)+p)-\vep/8+\vep_1\\
& < & f^{(m_k+1)}(y)-\vep/8+2\vep_1\leq f^{(m_k+1)}(y)
\end{eqnarray*}
and
\begin{eqnarray*}
s'+k\gamma+p & = &
(s+k\gamma)+(s'-s)+p>f^{(m_k)}(x)+p+\vep/8-\delta\\
& = & f^{(m_k)}(y)+(f^{(m_k)}(x)-f^{(m_k)}(y)+p)+\vep/8-\vep_1\\
& > & f^{(m_k)}(y)+\vep/8-2\vep_1\geq f^{(m_k)}(y).
\end{eqnarray*}
Thus
\[T^f_{k\gamma}(x,s)=(T^{m_k}x,s+k\gamma-f^{(m_k)}(x))\]
and
\[T^f_{k\gamma+p}(y,s')=(T^{m_k}y,s'+k\gamma+p-f^{(m_k)}(y)).\]
Hence
\begin{eqnarray*}
\lefteqn{d(T^f_{k\gamma}(x,s),T^f_{k\gamma+p}(y,s'))}\\ & =&
\|y-x\|+|(s'-s)+(f^{(m_k)}(x)-f^{(m_k)}(y)+p)|<
\delta+\vep_1<2\vep_1<\vep.
\end{eqnarray*}
It follows that
\begin{equation}\label{impli}
B\subset\{k\in\Z\cap[M,M+L):d(T^f_{k\gamma}(x,s),T^f_{k\gamma+p}(y,s'))<\vep\}.
\end{equation}
If $k\in(\Z\cap[M,M+L])\setminus B$ then either
$T^f_{k\gamma}(x,s)\notin X(\vep)$ or $m_k\notin A''$. Since for
every $m\in\N$ the set $\{k\in\N:m_k=m\}$ has at most $C/\gamma+1$
elements, we have
\[L-\#B \leq \#\left\{M\leq k<M+L:T^f_{k\gamma}(x,s)\notin
X(\vep)\right\}+\left(\frac{C}{\gamma}+1\right)(L'-\#A'').\] Hence
by (\ref{geba}), (\ref{nierapr}) and (\ref{nierl}) we obtain
\[L-\#B \leq
\frac{\vep}{2}L+\left(\frac{C}{\gamma}+1\right)4\vep_1L'\leq\left(\frac{\vep}{2}+4\frac{C+\gamma}{c}\vep_1\right)L\leq
\left(\frac{\vep}{2}+\frac{\vep}{2}\right)L\leq\vep L.\]
Consequently $(\#B)/L>1-\vep$, and (\ref{impli}) completes the
proof.
\end{proof}

\section{Ratner property for the von Neumann class of functions}
We call a function $f:\T\to\R$ {\em piecewise absolutely
continuous} if there exist $\beta_1,\ldots,\beta_k\in\T$ such that
$f|_{(\beta_j,\beta_{j+1})}$ is an absolutely continuous function
for $j=1,\ldots,k$ ($\beta_{k+1}=\beta_1$).
 Let
$d_j:=f_-(\beta_j)-f_+(\beta_j)$, where
$f_{\pm}(\beta)=\lim_{y\to\beta^\pm}f(y)$. Then the number
\begin{equation}\label{sumskok}
S(f):=\sum_{j=1}^kd_j=\int_{\T}f'(x)dx
\end{equation}
is the {\em sum of jumps} of $f$.

Let $T:\T\to\T$ be the rotation by an irrational number $\alpha$
which has bounded partial quotients. We will prove that if $f$ is
a positive piecewise absolutely continuous function with a
non--zero sum of jumps, then the special flow $T^f$ satisfies the
$\rat(t_0,P)$--property for every $t_0\neq 0$, where $P\subset
\R\setminus\{0\}$ is a non--empty finite set.

\begin{Lemma}\label{zababs}
Let $T:\T\to\T$ be the rotation by an irrational number $\alpha$
which has bounded partial quotients and let $f:\T\to\R$ be an
absolutely continuous function. Then
\[\sup_{0\leq n\leq q_{s+1}}\sup_{\|y-x\|<1/q_{s}}|f^{(n)}(y)-f^{(n)}(x)|\to 0\]
as $s\to\infty$.
\end{Lemma}

\begin{proof}
We first prove that if $T$ is an irrational rotation by $\alpha$
then
\begin{equation}\label{zbjed}
\sup_{0\leq n\leq
q_{s}}\sup_{\|y-x\|<1/q_{s}}|f^{(n)}(y)-f^{(n)}(x)|\to 0\;\;\text{
as }\;\;s\to\infty
\end{equation} for every
absolutely continuous $f:\T\to\R$. We recall that (\ref{zbjed})
was already proved to hold in \cite{Co-Fo-Si} (see Lemma~2 Ch.16,
\S3) for $C^1$--functions.

Let $f:\T\to\R$ be an absolutely continuous function. Then for
every $\vep>0$ there exists a $C^1$--function $f_{\vep}:\T\to\R$
such that
\[\sup_{x\in\T}|f(x)-f_{\vep}(x)|+\var(f-f_{\vep})<\vep/2.\]
Suppose that $0\leq n\leq q_{s}$ and $0<y-x<1/q_{s}$. Let us
consider the family of intervals
$\mathcal{I}=\{[x,y],[Tx,Ty],\ldots,[T^{n-1}x,T^{n-1}y]\}$. For
every $0\leq i\neq j<n$ we have
\[\|T^ix-T^jx\|\geq\|q_{s-1}\alpha\|>\frac{1}{2q_{s}},\]
by (\ref{ulla}). It follows that a point from $\T$ belongs to at
most two intervals from the family $\mathcal{I}$. Therefore
\begin{eqnarray*}\lefteqn{|(f^{(n)}(y)-f^{(n)}(x))-(f_{\vep}^{(n)}(y)-f_{\vep}^{(n)}(x))|}\\
&\leq&\sum_{i=0}^{n-1}|(f-f_{\vep})(T^iy)-(f-f_{\vep})(T^ix)|\\
&\leq&\sum_{i=0}^{n-1}\var_{[T^ix,T^iy]}(f-f_{\vep})\leq 2\var
(f-f_{\vep})<\vep.
\end{eqnarray*}
Since this inequality holds for every $\vep>0$ and the convergence
in (\ref{zbjed}) holds for $f_{\vep}$, by a standard argument, we
obtain (\ref{zbjed}) for $f$.

Suppose that $\alpha$ has bounded partial quotients and let
$C=\sup\{a_n:n\in\N\}+1$. Since every $0\leq n\leq q_{s+1}$ can be
represented as $n=bq_{s}+d$, where $b\leq a_{s+1}$ and $d\leq
q_{s-1}$, we have
\[\sup_{0\leq n\leq q_{s+1}}\sup_{\|y-x\|<1/q_{s}}|f^{(n)}(y)-f^{(n)}(x)|\leq
C\sup_{0\leq n\leq q_{s}}\sup_{\|y-x\|<1/q_{s}}|f^{(n)}(y)-f^{(n)}(x)|,\]
 which completes the proof.
\end{proof}

Let $T:\T\to\T$ be the rotation by an irrational number $\alpha$
which has bounded partial quotients and let
$C=\sup\{a_n:n\in\N\}+1$. Suppose that $f:\T\to\R$ is a positive
piecewise absolutely continuous function with a non--zero sum of
jumps $S=S(f)$. Put
\[D:=\{n_1d_1+\ldots+n_kd_k:0\leq n_1,\ldots,n_k\leq 2C+1\}.\]
 Since $D$ is
finite, we can choose $p\in(0,|S|)\setminus(D\cup(-D))$. Then
$0\notin\sgn(S)p-D$.
\begin{Th}\label{radlniez}
Suppose that $T:\T\to\T$ is the rotation by an irrational number
$\alpha$ with bounded partial quotients  and $f:\T\to\R$ a
positive and bounded away from zero piecewise absolutely
continuous function with a non--zero sum of jumps. Then the
special flow $T^f$ has the property $\rat(\gamma,(\sgn(S)p-D)\cup
(-\sgn(S)p+D))$ for every $\gamma>0$.
 \end{Th}

\begin{proof}
Without loss of generality we can assume that $f$ is continuous
from the right. A consequence of (\ref{sumskok}) is that we can
represent $f$ as the sum of two functions $f_{pl}$ and $f_{ac}$,
where $f_{ac}:\T\to\R$ is an absolutely continuous function with
zero mean and $f_{pl}:\T\to\R$ is piecewise linear with
$f_{pl}'(x)=S$ for all
$x\in\T\setminus\{\beta_1,\ldots,\beta_{k}\}$. The discontinuity
points and size of jumps of $f$ and $f_{pl}$ are the same.
Explicitly,
\[f_{pl}(x)=\sum_{i=1}^kd_i\{x-\beta_i\}+c\]
for some $c\in\R$.

 Let $C=\sup\{a_n:n\in\N\}+1$. Fix $\vep>0$ and $N\in\N$.
Then put
\[\kappa(\vep)=\frac{1}{k(2C+1)}\cdot\min\left(\frac{\vep}{2pC},\frac{1}{C^2}\right).\]
By Lemma~\ref{zababs}, there exists $s_0$ such that for any $s\geq
s_0$ we have
\begin{equation}\label{jnier0}
\sup_{0\leq n\leq
q_{s+1}}\sup_{\|y-x\|<1/q_{s}}|f_{ac}^{(n)}(y)-f_{ac}^{(n)}(x)|<\frac{\vep}{2}
\end{equation}
and
\begin{equation}\label{zalnan}
 \min\left(\kappa(\vep),1\right)\cdot
q_{s_0}>N.
\end{equation}
Then put
\[\delta(\vep,N)=\frac{p}{|S|q_{s_0+1}}.\]
 Take $x,y\in \T$ such that $0<\|x-y\|<\delta(\vep,N)$. Let $s$ be a (unique) natural number such that
\begin{equation}\label{dwapr}
\frac{p}{|S|q_{s+1}}<\|x-y\|\leq\frac{p}{|S|q_s}.
\end{equation}
Then $s\geq s_0$.
 Without loss of generality we can assume that $x<y$. We will also assume that $S>0$, in the case
$S<0$ the proof is similar. Let us consider the sequence
$\left(f_{pl}^{(n)}(y)-f_{pl}^{(n)}(x)\right)_{n\in\N}$. We have
\begin{eqnarray*}
\lefteqn{f_{pl}^{(n+1)}(y)-f_{pl}^{(n+1)}(x)}\\
&=&f_{pl}^{(n)}(y)-f_{pl}^{(n)}(x)+\sum_{i=1}^kd_i(\{y+n\alpha-\beta_i\}-\{x+n\alpha-\beta_i\})\\
&=&f_{pl}^{(n)}(y)-f_{pl}^{(n)}(x)+\sum_{i=1}^kd_i(y-x-\chi_{(x,y]}(\{\beta_i-n\alpha\})).
\end{eqnarray*}
It follows that  for every $n\geq 0$ we have
\begin{equation}\label{trow}
f_{pl}^{(n)}(y)-f_{pl}^{(n)}(x)=nS(y-x)-\overline{d}_n,
\end{equation}
where
\[\overline{d}_n:=\overline{d}_n(x,y)=\sum_{\{1\leq i\leq k,0\leq j<n:\{\beta_i-j\alpha\}\in(x,y]\}}d_i.\]
Take $1\leq i\leq k$. Suppose that
$\{\beta_i-k\alpha\},\{\beta_i-l\alpha\}\in(x,y]$, where $0\leq
k,l<q_{s+1}$ and $k\neq l$. Then
\[\|\{\beta_i-k\alpha\}-\{\beta_i-l\alpha\}\|\geq \|q_{s}\alpha\|>\frac{1}{2q_{s+1}}\geq\frac{1}{2Cq_{s}}.\]
It follows that the number of  discontinuities of
$f_{pl}^{(q_{s+1})}$ which are of the form $\beta_i-j\alpha$ and
are in the interval $(x,y]$ is less than
\[2Cq_s|y-x|+1\leq 2C\frac{p}{|S|}+1\leq 2C+1.\]
 It follows that the elements of the sequence
$(\overline{d}_n)_{n=1}^{q_{s+1}}$ belong to $D$. In view of
(\ref{trow}) and  (\ref{dwapr}) we have
\[f_{pl}^{(q_s)}(y)-f_{pl}^{(q_s)}(x)+\overline{d}_{q_s}=q_s S(y-x)\leq p\]
and
\[f_{pl}^{(q_{s+1})}(y)-f_{pl}^{(q_{s+1})}(x)+\overline{d}_{q_{s+1}}=q_{s+1}S(y-x)> p.\]
Moreover,  for any natural $n$ we have
\[0<f_{pl}^{(n+1)}(x)-f_{pl}^{(n+1)}(y)+\overline{d}_{n+1}-(f_{pl}^{(n)}(x)-f_{pl}^{(n)}(y)+\overline{d}_{n})
=S(y-x)\leq\frac{p}{q_{s}}.\]
 Hence,
there exists an integer interval $I\subset[q_s,q_{s+1}]$ such that
\[|f_{pl}^{(n)}(x)-f_{pl}^{(n)}(y)+\overline{d}_{n}-p|<\frac{\vep}{2}\text{ for }n\in  I\]
and
\[|I|\geq \min\left(\frac{\vep
}{2p}\,q_s,q_{s+1}-q_s\right)\geq
\min\left(\frac{\vep}{2pC},\frac{1}{C^2}\right)\cdot q_{s+1}.\]
Since $s\geq s_0$, by (\ref{jnier0}) we have
\[|f^{(n)}(x)-f^{(n)}(y)+\overline{d}_{n}-p|<\vep\text{ for }n\in  I.\]
Note that if $f^{(n)}$ and $f^{(n+1)}$ have the same points of
discontinuity in the interval $(x,y]$ then
$\overline{d}_n=\overline{d}_{n+1}$ and since $f^{(q_{s+1})}$ has
at most $k(2C+1)$ discontinuities in $(x,y]$, we can split $I$
into at most $k(2C+1)$ integer intervals on which the sequence
$(\overline{d}_n)_{n\in I}$ is constant. Thus we can choose $d\in
D$  and an integer subinterval $J\subset I$ such that
$\overline{d}_n=d$ for $n\in J$ and
\[|J|\geq\frac{1}{k(2C+1)}\cdot \min\left(\frac{\vep }{2pC},\frac{1}{C^2}\right)\cdot q_{s+1}
=\kappa(\vep)\cdot q_{s+1}.\] Therefore
\[|f^{(n)}(x)-f^{(n)}(y)-(p-d)|<\vep\text{ for }n\in J.\]
Now let $M,L$ be  natural numbers such that $J=[M,M+L]$. Then
\[\frac{L}{M}\geq\frac{|J|}{q_{s+1}}\geq\kappa(\vep),\]
\[M\geq q_s\geq q_{s_0}>N\;\;\text{ and }\;\;L\geq |J|\geq \kappa(\vep)q_{s+1}\geq\kappa(\vep)q_{s_0}>N,\]
by (\ref{zalnan}). Since the special flow  $T^f$ is weakly mixing
(see Proposition~2 in \cite{Iw-Le-Ma}), the automorphism
$T^f_{\gamma}$ is ergodic for all $\gamma\neq 0$, and hence an
application of Lemma~\ref{lemfundsp} completes the proof.
\end{proof}

Since special flows built over irrational rotations on the circle
and under piecewise absolutely continuous roof functions with a
non--zero sum of jumps are weakly mixing (see \cite{Iw-Le-Ma}),
from Theorems~\ref{prpr} and \ref{radlniez} we obtain the
following.

\begin{Cor}
Suppose that $T:\T\to\T$ is the rotation by an irrational number
$\alpha$ with bounded partial quotients and $f:\T\to\R$ is a
positive and bounded away from zero piecewise absolutely
continuous function with a non--zero sum of jumps. Then any
ergodic joining $\rho$ of the special flow $(T^f_t)_{t\in\R}$ and
an ergodic flow $(T_t)_{t\in\R}$ acting on $\ycn$ is either the
product joining, or $\rho$ is a finite extension of $\nu$.
\end{Cor}

\begin{problem}
It would be interesting to decide whether in the family of special
flows over the rotation by an irrational $\alpha$ with bounded
partial quotients and under $f$ which is piecewise absolutely
continuous with a non--zero sum of jumps we can find some with the
minimal self--joining property.
\end{problem}

\section{Absence of partial rigidity}\label{partrig}
\begin{Lemma}\label{lemsztyw}
Let $T:(X,\mathcal{B},\mu)\to(X,\mathcal{B},\mu)$ be an ergodic
automorphism and $f\in L^1(X,\mu)$ be a positive function such
that $f\geq c>0$. Suppose that the special flow $T^f$ is partially
rigid along a sequence $(t_n)$, $t_n\to+\infty$. Then there exists
$0<u\leq 1$ such that for every $0<\vep<c$ we have
\[\liminf_{n\to\infty}\mu\{x\in
X:\exists_{j\in\N}\;|f^{(j)}(x)-t_n|<\vep\}\geq u.\]
\end{Lemma}

\begin{proof}
By assumption, there exists $0<u\leq 1$ such that for any
measurable set $D\subset X^f$ we have
\[\liminf_{n\to\infty}\mu^f(D\cap T^f_{-t_n}D)\geq u\mu^f(D).\] Fix
$0<\vep<c$. Let ${A}:=X\times[0,\vep)$ and for any natural $n$ put
\[B_n=\{x\in
X:\exists_{j\in\N}\;|f^{(j)}(x)-t_n|<\vep\}.\]
Suppose that
$(x,s)\in{A}\cap T^f_{-t_n}{A}$. Then $0\leq s<\vep$ and there
exists $j\in\Z$ such that $0\leq s+t_n-f^{(j)}(x)<\vep$. It
follows that $-\vep<t_n-f^{(j)}(x)<\vep$, hence $x\in B_n$.
Therefore ${A}\cap T^f_{-t_n}{A}\subset B_n\times [0,\vep)$ and
\[\vep\liminf_{n\to\infty}\mu(B_n)=\liminf_{n\to\infty}\mu^f(B_n\times[0,\vep))\geq
\liminf_{n\to\infty}\mu^f({A}\cap T^f_{-t_n}{A})\geq
u\mu^f({A})=\vep u\] and the proof is complete.
\end{proof}

\begin{Th}\label{brakszt}
Let $T:\T\to\T$ be an ergodic rotation. Suppose that $f:\T\to\R$
is a positive and bounded away from zero piecewise absolutely
continuous function with $S(f)\neq 0$. Then the special flow $T^f$
is not partially rigid.
\end{Th}

\begin{proof}
Let $c$, $C$ be positive numbers such that $0<c\leq f(x)\leq C$
for every $x\in\T$. Let $\mu$ stand for Lebesgue measure on $\T$.
Assume, contrary to our claim, that $(t_n)$, $t_n\to+\infty$, is a
partial rigidity time for $T^f$. By Lemma~\ref{lemsztyw}, there
exists $0<u\leq 1$ such that for every $0<\vep<c$ we have
\begin{equation}\label{lemsztywwz}
\liminf_{n\to\infty}\mu\{x\in
\T:\exists_{j\in\N}\;|f^{(j)}(x)-t_n|<\vep\}\geq u.
\end{equation}
Without loss of generality we can assume that $S:=S(f)>0$, in the
case $S<0$ the proof is the same. Suppose that $\beta_i$,
$i=1,\ldots,k$ are all  points of discontinuity of $f$. Fix
\begin{equation}\label{okeps}
0<\vep<\min\left(\frac{Sc^2}{32kC(c+\var
f)+Sc^2)}u,\frac{c}{4}\right).
\end{equation}
Since $f'\in L^1(\T,\mu)$, there exists $0<\delta<\vep$ such that
$\mu(A)<\delta$ implies $\int_A|f'|\,d\mu<\vep$. Moreover, by the
ergodicity of $T$ (and recalling that $S=\int f'\,d\mu$) there
exist $A_\vep\subset\T$ with $\mu(A_\vep)>1-\delta$ and $m_0\in\N$
such that
\begin{equation}\label{szasr}
\frac{S}{2}\leq\frac{1}{m}f'^{(m)}(x)
\end{equation} for
all $m\geq m_0$ and $x\in A_\vep$.

Then take any $n\in\N$ such that $t_n/(2C)\geq m_0$ and
$t_n>2\vep$. Now
 let us consider the set $J_{n,\vep}$ of all
natural $j$ such that $|f^{(j)}(x)-t_n|<\vep$ for some $x\in\T$.
Then for such $j$ and $x$ we have
\[t_n+\vep>f^{(j)}(x)\geq
cj\;\;\text{ and } t_n-\vep<f^{(j)}(x)\leq Cj,\] whence
\begin{equation}\label{szaj}
t_n/(2C)\leq(t_n-\vep)/C<j<(t_n+\vep)/c\leq 2t_n/c
\end{equation}
for any $j\in J_{n,\vep}$; in particular,  $j\in J_{n,\vep}$
implies $j\geq m_0$.

Let $j_n=\max J_{n,\vep}$. The points of discontinuity of
$f^{(j_n)}$, i.e.\ $\{\beta_i-j\alpha\},1\leq i\leq k,0\leq
j<j_n$, divide  $\T$ into subintervals
$I_1^{(n)},\ldots,I_{kj_n}^{(n)}$. Some of these intervals can be
empty. Notice that for every $j\in J_{n,\vep}$ the function
$f^{(j)}$ is absolutely continuous on the interior of any interval
$I_i^{(n)}$, $i=1,\ldots,kj_n$.

Fix $1\leq i\leq kj_n$. For every $j\in J_{n,\vep}$ let
$I^{(n)}_{i,j}$ stand for the minimal closed subinterval of
$\overline{I^{(n)}_i}$ which includes the set $\{x\in
I^{(n)}_i:|f^{(j)}(x)-t_n|<\vep\}$. Of course, $I^{(n)}_{i,j}$ may
be empty. If $I^{(n)}_{i,j}=[z_1,z_2]$ is not empty then
\begin{equation}\label{malesza}
\left|\int_{I^{(n)}_{i,j}}\frac{f'^{(j)}}{j}\,d\mu\right|=\frac{|(f^{(j)})_-(z_2)-(f^{(j)})_+(z_1)|}{j}\leq\frac{2\vep}{j}
\leq\frac{4C\vep}{t_n}.
\end{equation}

Now suppose that $x$ is an end of $I_{i,j}^{(n)}$ and $y$ is an
end of $I_{i,j'}^{(n)}$ with $j\neq j'$. From (\ref{okeps}) it
follows that
\begin{equation}\label{wieksza}
\begin{aligned}
\int_{x}^y|f'|^{(j_n)}\,d\mu\geq &\left|\int_{x}^yf'^{(j)}\,d\mu\right|=|f^{(j)}(y)-f^{(j)}(x)|\\
\geq &|f^{(j)}(y)-f^{(j')}(y)|-|f^{(j')}(y)-t_n|-|f^{(j)}(x)-t_n|\\
\geq& c-2\vep\geq \frac{c}{2}.
\end{aligned}
\end{equation}

Let $K_i=\{j\in J_{n,\vep}:I_{i,j}^{(n)}\neq\emptyset\}$ and
suppose that $s=\#K_{i}\geq 1$. Then there exist $s-1$ pairwise
disjoint subintervals $H_{l}\subset I_{i}^{(n)}$, $l=1,\ldots,s-1$
that are disjoint from intervals $I_{i,j}^{(n)}$, $j\in K_i$ and
fill up the space between those intervals. In view of
(\ref{wieksza}) we have $\int_{H_l}|f'|^{(j_n)}\,d\mu\geq c/2$ for
$l=1,\ldots,s-1$. Therefore, by (\ref{malesza}) and
(\ref{wieksza}), we obtain
\begin{eqnarray}\label{konc1}
\begin{aligned}
\left|\sum_{j\in K_i}
\int_{I^{(n)}_{i,j}}\frac{f'^{(j)}}{j}\,d\mu\right|\;\;\;&\leq
&s\frac{4C\vep}{t_n}=\frac{4C\vep}{t_n}+\frac{8C\vep}{ct_n}(s-1)\frac{c}{2}\,\\
&\leq
&\frac{4C\vep}{t_n}+\frac{8C\vep}{ct_n}\sum_{l=1}^{s-1}\int_{H_l}|f'|^{(j_n)}\,d\mu\\
&\leq
&\frac{4C\vep}{t_n}+\frac{8C\vep}{ct_n}\int_{I^{(n)}_{i}}|f'|^{(j_n)}\,d\mu.\;\;\;
\end{aligned}
\end{eqnarray}
Since $\mu(A_{\vep}^c)<\delta$, we have
\begin{eqnarray}\label{konc2}
\begin{aligned}
\left|\sum_{i=1}^{kj_n}\sum_{j\in K_i}\int_{I^{(n)}_{i,j}\cap
A_\vep^c}\frac{f'^{(j)}}{j}\,d\mu\right|\;&\leq
&\frac{2C}{t_n}\sum_{i=1}^{kj_n}\sum_{j\in
K_i}\int_{I^{(n)}_{i,j}\cap A_\vep^c}|f'|^{(j_n)}\,d\mu\;\;\;\;\;\;\;\;\\
&\leq &\frac{2C}{t_n}\int_{A_\vep^c}|f'|^{(j_n)}\,d\mu\leq
\frac{2C}{t_n}j_n\vep\leq\frac{4C}{c}\vep.
\end{aligned}
\end{eqnarray}
As
\[B_n:=\{x\in\T:\exists_{j\in\N}\;|f^{(j)}(x)-t_n|<\vep\}\subset \bigcup_{i=1}^{kj_n}\bigcup_{j\in K_i
}I^{(n)}_{i,j},\] by (\ref{szasr}), (\ref{konc1}), (\ref{konc2})
and (\ref{szaj}) we have
\begin{eqnarray*}\frac{S}{2}\mu(B_n\cap A_\vep)&\leq &\sum_{i=1}^{kj_n}\sum_{j\in
K_i}\int_{I^{(n)}_{i,j}\cap A_\vep}\frac{f'^{(j)}}{j}\,d\mu\\
&\leq & \left|\sum_{i=1}^{kj_n}\sum_{j\in
K_i}\int_{I^{(n)}_{i,j}}\frac{f'^{(j)}}{j}\,d\mu\right|+\left|\sum_{i=1}^{kj_n}\sum_{j\in
K_i}\int_{I^{(n)}_{i,j}\cap
A_\vep^c}\frac{f'^{(j)}}{j}\,d\mu\right|\\&\leq &
kj_n\frac{4C\vep}{t_n}+\frac{8C\vep}{t_nc}\int_{\T}|f'|^{(j_n)}\,d\mu+\frac{4C}{c}\vep\\
&\leq &
\frac{8kC\vep}{c}+\frac{4C\vep}{c}+\frac{16C\vep}{c^2}\|f'\|_{L^1}\leq\frac{16kC}{c^2}(c+\var
f)\vep.
\end{eqnarray*}
Finally, from (\ref{okeps}) we obtain
\[\mu(B_n)\leq \mu(B_n\cap A_{\vep})+\mu(A_\vep^c)<\frac{32kC}{Sc^2}(c+\var
f)\vep+\vep<u,\] contrary to (\ref{lemsztywwz}).
\end{proof}

Collecting now Theorems~\ref{radlniez}, \ref{brakszt} and
Lemma~\ref{jmm} together with Remark~\ref{finext} we obtain the
following.

\begin{Th}\label{twmmsf}
Suppose that $T:\T\to\T$ is the rotation by an irrational number
$\alpha$ with bounded partial quotients and $f:\T\to\R$ is a
positive and bounded away from zero piecewise absolutely
continuous function with a non--zero sum of jumps. Then  the
special flow $(T^f_t)_{t\in\R}$ is mildly mixing.
\end{Th}

Applying  now the construction from Section~\ref{constr} we have
the following.

\begin{Cor}
On the two--dimensional torus there exists a $C^{\infty}$--flow
$(\varphi_t)_{t\in\R}$ with one singular point (of a simple pole
type) such that $(\varphi_t)$ preserves a positive
$C^{\infty}$--measure and the set of ergodic components of
$(\varphi_t)$ consists of a family of periodic orbits and one
non--trivial component of positive measure on which the flow is
mildly mixing but not  mixing. Moreover, on that component
$(\varphi_t)$ has the Ratner property $\rat(P)$ for some nonempty
finite set $P\subseteq\R\setminus\{0\}$.
\end{Cor}

\vspace{5mm}

\noindent Faculty of Mathematics and Computer Science\\ Nicolaus
Copernicus University\\ ul. Chopina 12/18\\ 87-100 Toru\'n, Poland

\vspace{3mm}

\noindent fraczek@mat.uni.torun.pl\\ mlem@mat.uni.torun.pl
\closegraphsfile

\end{document}